\RequirePackage{fix-cm}

\documentclass[smallextended]{svjour3}
\usepackage{doi}

\newtheorem{assumption}{Assumption}

\smartqed

\usepackage{amssymb, amsfonts, mathrsfs}
\usepackage{amsmath}
\usepackage{booktabs}

\usepackage{enumitem}
\setlist[enumerate,1]{label=\textnormal{(\roman*)}}

\usepackage[numbers]{natbib}

\newcommand{\cO}{\mathcal{O}}
\newcommand{\cL}{\mathscr{L}}
\newcommand{\cD}{\mathscr{D}}
\newcommand{\cM}{\mathcal{M}}
\newcommand{\cA}{\mathcal{A}}
\newcommand{\1}{\mathbf{1}}
\newcommand{\E}{\mathbb{E}}
\newcommand{\R}{\mathbb{R}}

\begin{document}

\title{Compatibility and Exclusivity of Loss Functions under Statistical Optimality}
\subtitle{\,}

\author{Stanisław M. S. Halkiewicz}

\institute{S. M. S. Halkiewicz \at
              Group of Machine Learning Research, Jagiellonian University, Kraków, Poland \\
              \email{stashal@o2.pl}
}

\date{Received: date / Accepted: date}

\maketitle

\begin{abstract}
In applications a loss is usually chosen from a family, like a quantile level, a robustness threshold or a surrogate margin loss. In this paper, we ask which of these choices change the optimal decision. A model, a class of decision rules and an
optimality criterion define an optimality map $L\mapsto\cO(L)$. Two losses are compatible when they share an optimal rule, and exclusivity regions are unions of the connected components of this relation. When $\cO$ is single-valued these components are the fibres of the map, and the problem is to identify the
statistically meaningful features of a loss that determine its fibre. Under explicit model assumptions, the answer is the quantile level for generalized piecewise linear quantile scores and the margin invariant (equivalently, the link) for regular margin losses read at the score level. For canonical Huber losses on strictly skewed location models, a class that includes every asymmetric Laplace model, it is the threshold. Two structural results hold in general. Convex interpolation between losses with distinct unique Bayes acts meets
uncountably many compatibility classes, so the structure is conic rather than convex. Coarsening the optimality operator can only merge components, which is how the known common Bayes classifier of calibrated surrogates arises.

\end{abstract}

\section{Introduction}
\label{intro}

Statistical procedures are chosen relative to a loss; in practice, however, the loss is rarely fixed. A forecaster might chooses a quantile level, a robust regression analyst chooses a Huber threshold, a classifier is trained on one of several
convex surrogates, and losses are rescaled, then reweighted and mixed. Each such choice moves the loss within a family. In this paper, we ask which of these moves change the optimal decision and which leave it unchanged. The answer depends on the statistical model, on the class of decision rules and on the optimality criterion. For a symmetric posterior, squared and
absolute error have the same Bayes act; for a skewed one they do not.

Several pieces of literature answer parts of this question. Classical decision theory studies optimality for a fixed loss \cite{Wald1950,Berger1985,LehmannCasella1998}.
Elicitation fixes a functional and asks which losses have it as their Bayes act. Quantiles are elicited exactly by the generalized piecewise linear scores \cite{Gneiting2011Quantiles,Gneiting2011}, each of which is a mixture of elementary scores \cite{EhmEtAl2016}. Osband's principle relates consistent scores to identification functions
\cite{Osband1985,GneitingRaftery2007,FisslerZiegel2016,SteinwartEtAl2014,Dimitriadis2024}, and the first-order conditions in our proofs are identification functions in this (Osband's) sense. In classification, calibration theory shows that convex surrogates share the Bayes classifier \cite{Zhang2004,BartlettJordanMcAuliffe2006}, and
proper composite loss theory describes how a Bayes score encodes the class probability through a link \cite{ReidWilliamson2010}. Properization \cite{BrehmerGneiting2020} composes a scoring rule with the Bayes-act map, which constitutes our basic object. M-estimation \cite{Huber1964,HampelEtAl1986} and M-quantiles \cite{BrecklingChambers1988} supply the robust families.

These results are organized around a target: a functional, a classifier or a probability. They do not, by themselves, organize a \emph{given} family of losses into regimes of mutually compatible optima relative to an arbitrary operator and
decision class. Nor do they ask how such an organization behaves under scaling, interpolation or coarsening of the operator, and these are the questions taken up here. Most of the ingredients are classical; what is new are the organizing statements and a few results that the target-centred view does not naturally produce.

A model, a decision class $\cD$ and an optimality
criterion define an optimality map $L\mapsto\cO(L)\subseteq\cD$. Two losses are \emph{compatible} if $cO(L)\cap\cO(L')\neq\emptyset$. A set of losses is an \emph{exclusivity region} if no loss inside it is compatible with one outside it.
When $\cO$ is single-valued, the natural classes are just the fibres of $L\mapsto\cO(L)$ (Lemma~\ref{lem:invariant},
Proposition~\ref{prop:components}). The results are of two kinds: loss-side invariants that describe these fibres exactly under explicit model assumptions, and statements on how the resulting structure behaves under natural operations on losses and operators.

Section~\ref{sec:families} identifies such invariants for three families. Under posterior Bayes optimality, the quantile level is a complete invariant on generalized piecewise linear quantile scores (Theorem~\ref{thm:quantile}). The margin invariant, equivalently the link, is a complete invariant for regular margin losses at the score level when the class probability has full support (Theorem~\ref{thm:margin}). For canonical Huber losses, distinct thresholds have distinct optimal equivariant estimators on strictly skewed location models (Theorem~\ref{thm:huber}). This model class contains every asymmetric Laplace model (Proposition~\ref{prop:laplace}), and nothing is claimed for other losses sharing a quadratic-core radius. Section~\ref{sec:structure} treats the general structure: exclusivity regions are unions of compatibility components (Proposition~\ref{prop:components}); Bayes-type operators ignore scaling and shifts in $\theta$ Proposition~\ref{prop:invariance}, essentially standard); convex interpolation between losses with distinct unique Bayes acts meets uncountably many compatibility classes (Theorem~\ref{thm:convexity}); and coarsening the operator can only merge components (Proposition~\ref{prop:coarsening}).

Several of the facts used are known. Generalized piecewise linear scores are strictly consistent for quantiles \cite{Gneiting2011}, and Theorem~\ref{thm:quantile} restates this as a complete-invariant statement without adding to quantile elicitation itself. Each margin loss has a Bayes score given by its link \cite{Zhang2004,ReidWilliamson2010}, calibrated surrogates share the Bayes classifier \cite{BartlettJordanMcAuliffe2006}, and positive scaling does not change optimality. As far as we are aware, the following are new:
\begin{itemize}
\item the treatment of compatibility through components relative to a triple (loss family, decision class, operator) (Proposition~\ref{prop:components});
\item the convex interpolation theorem (Theorem~\ref{thm:convexity});
\item the coarsening principle and the readout collapses it organizes
(Proposition~\ref{prop:coarsening}, Section~\ref{sec:collapse});
\item the complete-invariant form of Theorem~\ref{thm:margin}, including the
role of richness;
\item the injectivity of $\kappa\mapsto b_\kappa$ for canonical Huber estimators under strict skewness, with an analytic verification for asymmetric Laplace noise (Theorem~\ref{thm:huber}, Proposition~\ref{prop:laplace}).
\end{itemize}
The mere existence of trivial exclusivity partitions is not among the contributions, since such partitions always exist (Proposition~\ref{prop:components}). Instead, the interest lies in those that are described by meaningful invariants.

After the setup in Section~\ref{sec:setup}, Sections~\ref{sec:families}
and~\ref{sec:structure} contain the results above, Section~\ref{sec:collapse}
gives the collapse examples, and Section~\ref{sec:discussion} discusses
limitations.

\section{Optimality maps and exclusivity}
\label{sec:setup}

\paragraph{Notation.}
Let $\{P_\theta:\theta\in\Theta\}$ be a statistical model on a sample space $(\mathcal{X},\mathcal{F})$, with $\Theta\subseteq\R^d$ and densities $f(x;\theta)$ jointly measurable with respect to a $\sigma$-finite dominating measure. The action space $\cA$ is a measurable subset of a finite-dimensional normed space; in Section~\ref{sec:families}, $\cA=\R$. A decision rule is a measurable $\delta:\mathcal{X}\to\cA$, and $\cD$ denotes the class of rules under consideration. A loss is a measurable $L:\Theta\times\cA\to[0,\infty)$ with $a\mapsto L(\theta,a)$ continuous for every $\theta$; $\cL$ denotes the set of all losses. The risk is $R_L(\theta,\delta)=\E_\theta L(\theta,\delta(X))\in[0,\infty]$.
For Bayesian statements we fix a prior $\Pi$, write $P_X$ for the marginal law of $X$ and $\Pi_x$ for a regular conditional distribution (posterior) of $\theta$ given $X=x$, so that $(x,t)\mapsto\Pi_x(\theta<t)$ is jointly measurable. Rules that agree $P_X$-almost surely are identified.

\begin{definition}[Optimality map; compatibility]
\label{def:optimality}
An \emph{optimality operator} is a map $\cO:\cL\to2^{\cD}$; $\cO(L)$ is the (possibly empty) set of rules regarded as optimal for $L$. Two losses $L,L'$ are \emph{$\cO$-compatible} if $\cO(L)\cap\cO(L')\neq\emptyset$.
\end{definition}

Three operators are used below. \emph{Posterior Bayes optimality}
\begin{equation}
\label{eq:posterior-bayes}
\cO_\Pi(L):=\Bigl\{\delta\in\cD:\ \delta(x)\in\arg\min_{a\in\cA}\int L(\theta,a)\,\Pi_x(d\theta)
\ \text{for $P_X$-a.e.\ }x\Bigr\}
\end{equation}
coincides with minimization of the Bayes risk when the latter is finite, and remains meaningful when it is not; when $X$ carries no information, $\Pi_x=Q$ for a fixed law $Q$ and $\cO_Q(L)=\arg\min_{a}\int L(\theta,a)\,Q(d\theta)$ is the set of Bayes acts. \emph{Oracle optimality over a subclass} $\cD_0\subseteq\cD$ selects the $\delta\in\cD_0$ with $R_L(\theta,\delta)\le R_L(\theta,\delta')$ for all $\theta$ and all $\delta'\in\cD_0$. \emph{Minimax optimality} selects minimizers of $\sup_\theta R_L(\theta,\delta)$. Compatibility is symmetric, and reflexive on $\{L:\cO(L)\neq\emptyset\}$, but it need not be transitive when $\cO$ is set-valued (Section~\ref{subsec:components}).

\begin{definition}[Exclusivity]
\label{def:exclusivity}
Fix an ambient family $\cM\subseteq\cL$. A set $\mathcal{C}\subseteq\cM$ is an \emph{exclusivity region within $\cM$} if no loss in $\mathcal{C}$ is $\cO$-compatible with a loss in $\cM\setminus\mathcal{C}$. A family of pairwise disjoint exclusivity regions $\{\mathcal{C}_i\}_{i\in I}$ in $\cM$ is an \emph{exclusivity partition} of $\cM$; it is \emph{total} if $\bigcup_i\mathcal{C}_i=\cM$, and \emph{realizable} if every $\mathcal{C}_i$ contains some $L$ with $\cO(L)\neq\emptyset$.
\end{definition}

Realizability excludes a trivial case: a region of losses that have no optimal rule at all is formally exclusive but says nothing about optimal behaviour.
The single-class partition $\{\cM\}$ is always a total exclusivity partition and carries no information. Exclusivity is used below in the following form, which is immediate.

\begin{corollary}[No rule serves two regions]
\label{cor:nfl}
If $\{\mathcal{C}_i\}$ is an exclusivity partition of $\cM$, $L\in\mathcal{C}_i$ and $L'\in\mathcal{C}_j$ with $i\neq j$, then $\cO(L)\cap\cO(L')=\emptyset$.
\end{corollary}

\begin{proof}
$L'\in\cM\setminus\mathcal{C}_i$ by disjointness, and $\mathcal{C}_i$ is an exclusivity region within $\cM$. \qed
\end{proof}

\paragraph{Fibres and invariants.}
If $\cO$ is single-valued on $\cM$, compatibility is the equivalence relation $\cO(L)=\cO(L')$, and its classes are the fibres of $L\mapsto\cO(L)$. In that case the finest total exclusivity partition is simply the partition into fibres (Proposition~\ref{prop:components} treats the general case). This description is correct, but it says nothing about which features of $L$ determine the fibre.
We ask whether a natural loss-side quantity, such as an asymmetry level, a margin invariant or a threshold, determines the fibre, and under which model assumptions. Lemma~\ref{lem:invariant} reduces this to a separation property, and each theorem in Section~\ref{sec:families} is proved by verifying that property. On terminology: \emph{compatible} means sharing an optimal rule, and a \emph{compatibility component} is a connected component of the compatibility relation (Section~\ref{subsec:components}). We speak of \emph{compatibility classes} only when compatibility is an equivalence relation, in particular when $\cO$ is single-valued; the classes are then the fibres.

\begin{lemma}[Invariants and exclusivity]
\label{lem:invariant}
Let $\cM\subseteq\cL$, let $I:\cM\to\Lambda$ be any map, and put
$\cM_\lambda:=I^{-1}(\lambda)$ for $\lambda\in I(\cM)$.
\begin{enumerate}
\item $\{\cM_\lambda\}$ is a total exclusivity partition of $\cM$ if and only if
\[
I(L)\neq I(L')\ \Longrightarrow\ \cO(L)\cap\cO(L')=\emptyset
\qquad (L,L'\in\cM). \tag{S}
\]
It is realizable if and only if each $\cM_\lambda$ contains some $L$ with $\cO(L)\neq\emptyset$.
\item If (S) holds and $\cO$ is nonempty and constant on each $\cM_\lambda$, then the $\cM_\lambda$ are exactly the fibres of $L\mapsto\cO(L)$ on $\cM$. We then call $I$ a \emph{complete invariant}.
\end{enumerate}
\end{lemma}

\begin{proof}
(i) The level sets are disjoint and cover $\cM$. $\cM_\lambda$ is an exclusivity region within $\cM$ if and only if no $L\in\cM_\lambda$ is compatible with any $L'\in\cM$ with $I(L')\neq\lambda$; requiring this for every $\lambda$ is (S).
The realizability statement is the definition. (ii) $\cO$ is constant on each level set, and by (S) the nonempty values on distinct level sets are disjoint, hence distinct. \qed
\end{proof}

\section{Loss families and incompatible Bayes acts}
\label{sec:families}

Throughout this section $\cA=\R$ and the parameter (or label) is scalar. In each family, optimality reduces to one equation whose form is fixed by a single loss-side quantity, and Lemma~\ref{lem:invariant} turns the reduction into a partition statement.

\subsection{Quantile scores and the asymmetry level}
\label{subsec:quantile}

Quantiles are the Bayes acts of asymmetric linear loss \cite{KoenkerBassett1978,Koenker2005}.
More generally, under mild regularity the scoring functions consistent for the $\tau$-quantile are exactly the generalized piecewise linear scores $(\1\{\theta<a\}-\tau)(g(a)-g(\theta))$ with $g$ nondecreasing, and they are strictly consistent when $g$ is strictly increasing \cite{Gneiting2011Quantiles,Gneiting2011}. Ehm et al.\ \cite{EhmEtAl2016} represent each of them as a mixture, with mixing measure $dg$, of elementary scores. The invariant is the level $\tau$; the Bayes act depends neither on the free parameter $g$ nor on an additive term in $\theta$.

Let $\Theta\subseteq\R$ contain at least two points. For $\tau\in(0,1)$, a continuous strictly increasing $g:\R\to\R$ and a measurable $h:\Theta\to[0,\infty)$ put
\begin{equation}
\label{eq:gpl}
L_{\tau,g,h}(\theta,a):=(\1\{\theta<a\}-\tau)\bigl(g(a)-g(\theta)\bigr)+h(\theta),
\end{equation}
and let $\cL_{\mathrm q}$ be the set of all such losses. The check loss $\rho_\tau(\theta-a)$, $\rho_\tau(u)=(\tau-\1\{u<0\})u$, is $L_{\tau,\mathrm{id},0}$.
Each $L_{\tau,g,h}$ is nonnegative and continuous in $a$. The level is determined by the loss - if $L_{\tau,g,h}=L_{\tau',g',h'}$, then $h=h'$ (take $a=\theta$),
and for $\theta_0<\theta_1$ in $\Theta$ the pairs $(\theta_0,\theta_1)$ and
$(\theta_1,\theta_0)$ give
$g'(\theta_1)-g'(\theta_0)=\frac{1-\tau}{1-\tau'}(g(\theta_1)-g(\theta_0))=\frac{\tau}{\tau'}(g(\theta_1)-g(\theta_0))$, whence $\tau=\tau'$. We write $\tau(L)$ and $\cL_\tau:=\{L\in\cL_{\mathrm q}:\tau(L)=\tau\}$. In general $L_{\tau,g,h}$ is not convex in $a$. If $\Theta=\R$, convexity for every $\theta$ forces $g$ to be convex on every $(\theta,\infty)$ and concave on every $(-\infty,\theta)$, hence affine. A convexity requirement would therefore reduce $\cL_\tau$ to $\{c\rho_\tau(\theta-a)+h(\theta):c>0\}$, and the argument below does not need one.

\begin{assumption}[Posterior regularity]
\label{ass:posterior}
For $P_X$-a.e.\ $x$, the posterior $\Pi_x$ is supported on an interval $S_x$ and its distribution function $F_x$ is continuous on $\R$ and strictly increasing on $S_x$.
\end{assumption}

Continuity rules out atoms, and strict increase on an interval rules out gaps. Thanks to these, every level $\tau\in(0,1)$ have a unique posterior quantile $q_x(\tau):=F_x^{-1}(\tau)$ in the interior of $S_x$. Both parts are needed: a mixed posterior such as $\tfrac12\delta_0+\tfrac12\mathrm{Unif}[1,2]$ assigns the act $0$ to every $\tau<\tfrac12$: atoms make distinct levels share an act. Gaps (flat parts of $F_x$) make Bayes acts set-valued, and together with atoms they can destroy transitivity of compatibility (Section~\ref{subsec:components}).

\begin{theorem}[Quantile levels]
\label{thm:quantile}
Let $\cO_\Pi$ be posterior Bayes optimality~\eqref{eq:posterior-bayes} with
$\cA=\R$, and assume Assumption~\ref{ass:posterior}. Let $\cL_{\mathrm q}^\Pi$ be the set of $L_{\tau,g,h}\in\cL_{\mathrm q}$ with $\E[\,|g(\theta)|+h(\theta)\mid X=x]<\infty$ for $P_X$-a.e.\ $x$ (for the check loss, $\E[|\theta|\mid X]<\infty$ a.s.). Assume $\cD$ contains the quantile rules $\delta_\tau(x):=q_x(\tau)$, $\tau\in(0,1)$. Then $\cO_\Pi(L)=\{\delta_{\tau(L)}\}$ for every $L\in\cL_{\mathrm q}^\Pi$. Consequently $\tau$ is a complete invariant on $\cL_{\mathrm q}^\Pi$: the sets
$\cL_\tau\cap\cL_{\mathrm q}^\Pi$, $\tau\in(0,1)$, form a realizable total exclusivity partition of $\cL_{\mathrm q}^\Pi$, and they are exactly the fibres of $L\mapsto\cO_\Pi(L)$.
\end{theorem}

\begin{proof}
Fix $L=L_{\tau,g,h}\in\cL_{\mathrm q}^\Pi$ and $x$ outside a $P_X$-null set, and write $F=F_x$, $q=q_x(\tau)$ and $\varphi(a):=\E[L(\theta,a)\mid X=x]$, which is finite because $|L(\theta,a)|\le|g(a)|+|g(\theta)|+h(\theta)$. For $b<a$ and every
$\theta$,
\[
L(\theta,a)-L(\theta,b)=\int_{(b,a]}\bigl(\1\{\theta<t\}-\tau\bigr)\,dg(t),
\]
as one checks separately for $\theta<b$, $b\le\theta<a$ and $\theta\ge a$; here $dg$ is the Lebesgue--Stieltjes measure of $g$, which has no atoms and charges every nondegenerate interval. The integrand is bounded by $1$ and
$dg((b,a])<\infty$, so Fubini's theorem gives
\begin{equation}
\label{eq:quantile-increment}
\varphi(a)-\varphi(b)=\int_{(b,a]}\bigl(F(t)-\tau\bigr)\,dg(t),
\qquad F(t)=\Pi_x(\theta<t),
\end{equation}
using continuity of $F$. If $t<q$ then $F(t)<\tau$: either $t\in S_x$ and strict monotonicity applies, or $t$ lies below $S_x$ and $F(t)=0$. Similarly $F(t)>\tau$ for $t>q$. Hence, by~\eqref{eq:quantile-increment}, $\varphi(a)<\varphi(b)$ whenever $b<a\le q$, and $\varphi(a)>\varphi(b)$ whenever $q\le b<a$, so $q$ is the unique minimizer of $\varphi$ on $\R$. Thus every element of $\cO_\Pi(L)$ equals $\delta_\tau$ $P_X$-a.s., and $\delta_\tau\in\cD$ is measurable, being the generalized inverse of a jointly measurable distribution function. So $\cO_\Pi(L)=\{\delta_\tau\}$. If
$\tau\neq\tau'$ then $q_x(\tau)\neq q_x(\tau')$ for $P_X$-a.e.\ $x$, so
$\{\delta_\tau\}\cap\{\delta_{\tau'}\}=\emptyset$. This is (S) for $I=\tau$;
$\cO_\Pi$ is nonempty and constant on each level set, and
Lemma~\ref{lem:invariant} gives the conclusion. \qed
\end{proof}

\paragraph{Interpretation.}
Distinct levels correspond to distinct posterior quantiles, so losses of different levels cannot share a Bayes rule under the stated conditions. Within a level, the Bayes rule depends on no other feature of the loss. The compatibility class of the check loss $\rho_\tau$ consists of all strictly consistent quantile scores of level $\tau$, whose analytic forms can differ considerably (for instance, $g(t)=t$ and $g(t)=e^t$), together with arbitrary integrable additive terms $h(\theta)$. In decision-theoretic terms, this is strict consistency \cite{Gneiting2011} expressed as a partition of the loss family.

The theorem says nothing about \emph{quantile crossing}. The true quantile rules $\tau\mapsto\delta_\tau(x)$ are monotone and never cross.
Theorem~\ref{thm:quantile} explains only why several levels require several rules: no single rule is Bayes for two levels (Corollary~\ref{cor:nfl}). If these rules are estimated separately, the estimates may cross \cite{Koenker2005,BondellReichWang2010}, but this is a property of the estimation procedure and therefore is not a consequence of exclusivity.

\paragraph{Mixtures across levels.}
Let $\tau_0<\tau_1$. Since $\rho_\tau(u)$ is affine in $\tau$,
$(1-w)\rho_{\tau_0}+w\rho_{\tau_1}=\rho_{\tau_w}$ with $\tau_w=(1-w)\tau_0+w\tau_1$.
More generally $(1-w)L_{\tau_0,g,h_0}+wL_{\tau_1,g,h_1}=L_{\tau_w,g,(1-w)h_0+wh_1}$ for a common $g$. By Theorem~\ref{thm:quantile}, a weighted pinball objective does not serve both levels: its Bayes rule is the $\tau_w$-quantile, which for $w\in(0,1)$ differs a.s.\ from both endpoint quantiles. For a standard exponential predictive law and $(\tau_0,\tau_1)=(0.25,0.75)$, the equally weighted objective targets the median $0.693$, whereas the average of $q(0.25)=0.288$ and $q(0.75)=1.386$ is $0.837$. By Section~\ref{subsec:convexity}, this passage through a continuum of classes occurs for convex interpolation in general.

\paragraph{A distribution-dependent route.}
Gzyl \cite{Gzyl2024} reaches quantiles from a different direction. For a fixed law with continuous, strictly increasing distribution function $F$ and a suitable transformation $h$ of $[0,1]$, the best predictor in the transformed squared distance $\E\bigl(h(F(\theta))-h(F(a))\bigr)^2$ is the quantile of $F$ at a level determined by $h$. For that fixed distribution, this yields a further family of losses grouped by their induced quantile level. Since these losses are built from $F$ itself, the construction is not a distribution-free elicitation result, and we do not develop that model-dependent family here.

\subsection{Margin losses and the score level}
\label{subsec:margin}

In binary classification with label $Y\in\{-1,+1\}$ (playing the role of $\theta$) and score rules $f:\mathcal{X}\to\R$, a margin loss is $L_\phi(y,t)=\phi(yt)$. Classification calibration
\cite{Zhang2004,BartlettJordanMcAuliffe2006} concerns the sign of the Bayes score. Proper composite loss theory \cite{ReidWilliamson2010} concerns the link through which the Bayes score encodes $\eta(x):=\mathbb{P}(Y=+1\mid X=x)$. We work at the level of the score itself.

Call $\phi:\R\to[0,\infty)$ \emph{regular} if it is convex, continuously differentiable and strictly decreasing, and its \emph{margin invariant}
\[
H_\phi(t):=\frac{\phi'(-t)}{\phi'(t)}
\]
is a continuous strictly increasing bijection $\R\to(0,\infty)$. Here $\phi'<0$ everywhere, as if $\phi'(t_0)=0$, convexity would make $\phi$ nondecreasing on $[t_0,\infty)$. Let $\cL_{\mathrm m}$ be the set of regular margin losses and, for $\lambda$ in $\Lambda:=\{H_\phi:L_\phi\in\cL_{\mathrm m}\}$, let $\cL_\lambda:=\{L_\phi\in\cL_{\mathrm m}:H_\phi=\lambda\}$. The invariant is unchanged by positive scaling of $\phi$. Pointwise minimization of the conditional risk $C_\eta(t)=\eta\phi(t)+(1-\eta)\phi(-t)$ gives the stationarity condition $\eta\phi'(t)=(1-\eta)\phi'(-t)$ \cite{Zhang2004,BartlettJordanMcAuliffe2006},
that is,
\begin{equation}
\label{eq:margin-foc}
H_\phi\bigl(f^*(x)\bigr)=\frac{\eta(x)}{1-\eta(x)},\qquad \eta(x)\in(0,1).
\end{equation}
For $\eta\in(0,1)$, $C_\eta$ is convex and $C^1$, and by~\eqref{eq:margin-foc} it has the unique minimizer $H_\phi^{-1}(\eta/(1-\eta))$.

\begin{assumption}[Richness of the class probability]
\label{ass:richness}
The class probability satisfies $\eta(X)\in(0,1)$ $P_X$-a.s., and
$P_X(\eta(X)\in U)>0$ for every nonempty open interval $U\subseteq(0,1)$.
\end{assumption}

\begin{theorem}[Margin invariants at the score level]
\label{thm:margin}
Let $\cD$ be the set of measurable score rules and $\cO_{\mathrm{score}}$ posterior Bayes optimality, that is, $f\in\cO_{\mathrm{score}}(L_\phi)$ if and only if $f(x)$ minimizes $C_{\eta(x)}$ for $P_X$-a.e.\ $x$. Under Assumption~\ref{ass:richness}, $\cO_{\mathrm{score}}(L_\phi)=\{H_\phi^{-1}\circ(\eta/(1-\eta))\}$ for every $L_\phi\in\cL_{\mathrm m}$, and $H_\phi$ is a complete invariant, as the classes $\{\cL_\lambda\}_{\lambda\in\Lambda}$ form a realizable total exclusivity partition of $\cL_{\mathrm m}$, consisting of the fibres of $L\mapsto\cO_{\mathrm{score}}(L)$.
\end{theorem}

\begin{proof}
On the $P_X$-full event $\{\eta\in(0,1)\}$ the minimizer of $C_{\eta(x)}$ is unique and equal to $f^*_\phi(x)=H_\phi^{-1}(\eta(x)/(1-\eta(x)))$. This is measurable because $H_\phi^{-1}$ is continuous, so $\cO_{\mathrm{score}}(L_\phi)=\{f^*_\phi\}$, which is nonempty and constant on each $\cL_\lambda$. For (S), let $H:=H_\phi\neq H':=H_{\phi'}$ and suppose that $f$ is optimal for both. Then $H(f(x))=\eta(x)/(1-\eta(x))=H'(f(x))$ a.s. By continuity there is a nonempty open interval $V$ on which $H\neq H'$. The set $H(V)$ is a nonempty open subinterval of $(0,\infty)$, so by Assumption~\ref{ass:richness}
the event $\{\eta(X)/(1-\eta(X))\in H(V)\}$ has positive probability. On it, $f(X)=H^{-1}(\eta/(1-\eta))\in V$, hence $H(f(X))\neq H'(f(X))$, a contradiction. Lemma~\ref{lem:invariant} completes the proof. \qed
\end{proof}

\paragraph{Relation to links.}
In the language of \cite{ReidWilliamson2010}, $L_\phi$ is a composite loss with link $\psi_\phi(\eta)=H_\phi^{-1}(\eta/(1-\eta))$, and $H_\phi$ and $\psi_\phi$ determine each other. That each margin loss has such a link, and that the link gives the Bayes score, is known. Theorem~\ref{thm:margin} adds that, under richness, the link is a \emph{complete invariant} of score-level compatibility - losses with different links never share an optimal score rule, and losses with the same link always do. Richness cannot be dropped. If $\eta\equiv c$ is constant, the Bayes score is $f\equiv H_\phi^{-1}(c/(1-c))$, so any two invariants that agree at that one point share it; for instance, with $\eta\equiv\tfrac12$ every logistic loss $\phi_\lambda$ below has Bayes score $f\equiv0$.

\paragraph{The logistic family.}
For $\phi_\lambda(t)=\log(1+e^{-\lambda t})$, $\lambda>0$, one finds $H_{\phi_\lambda}(t)=e^{\lambda t}$, so $f^*_\lambda=\lambda^{-1}\,\mathrm{logit}\,\eta$. Distinct $\lambda$ give distinct score rules wherever $\eta\neq\tfrac12$, an event of positive probability under Assumption~\ref{ass:richness}. All $f^*_\lambda$ nevertheless have the same sign $\mathrm{sign}(\eta-\tfrac12)$, and each recovers the probability exactly, $\eta=(1+e^{-\lambda f^*_\lambda})^{-1}$. The family is therefore exclusive when the raw score is read, and collapsed when only the classifier or the recovered probability is read. In Section~\ref{subsec:margin-collapse} this contrast is explained by the coarsening principle.

\subsection{The canonical Huber family on strictly skewed location models}
\label{subsec:huber}

For $\kappa>0$ let $\ell_\kappa(u)=\tfrac12u^2$ for $|u|\le\kappa$ and $\ell_\kappa(u)=\kappa|u|-\tfrac12\kappa^2$ otherwise, with derivative $\psi_\kappa(u)=\max(-\kappa,\min(u,\kappa))$, and put $L_\kappa(\theta,a):=\ell_\kappa(a-\theta)$. In the location model $X=\theta+Z$, consider the equivariant rules $\delta_b(x)=x-b$, $b\in\R$, forming $\cD_{\mathrm{eq}}$. The risk of $\delta_b$ is $\E\,\ell_\kappa(Z-b)$ for every $\theta$, so oracle optimality $\cO_{\mathrm{eq}}$ over $\cD_{\mathrm{eq}}$ selects the minimizing offsets. M-estimation theory \cite{Huber1964,HampelEtAl1986,vanDerVaart1998}
characterizes these by the estimating equation
\begin{equation}
\label{eq:huber-foc}
g(b,\kappa):=\E\,\psi_\kappa(Z-b)=0.
\end{equation}
The offset $b_\kappa$ defines a location functional that interpolates between the median ($\kappa\to0$) and the mean ($\kappa\to\infty$). We ask whether distinct thresholds define distinct functionals. If $Z$ is symmetric about $b_0$, then $g(b_0,\kappa)=0$ for every $\kappa$ by oddness of $\psi_\kappa$, and all thresholds collapse. Some asymmetry is therefore necessary.

\begin{definition}[Strictly skewed noise]
\label{def:skewed}
Let $Z$ have a continuous, strictly positive density and $\E|Z|<\infty$. Then $\partial_bg(b,\kappa)=-P(|Z-b|\le\kappa)<0$, and $g(\cdot,\kappa)$ decreases from $\kappa$ to $-\kappa$, so~\eqref{eq:huber-foc} has a unique root $b_\kappa$.
$Z$ is \emph{strictly skewed} if there is $\sigma\in\{+1,-1\}$ such that $\sigma\bigl(P(Z-b_\kappa>\kappa)-P(Z-b_\kappa<-\kappa)\bigr)>0$ for every $\kappa>0$.
\end{definition}

The condition is imposed along the root curve $b=b_\kappa$. A version \textit{for every $b$} would be unsatisfiable, since for a full-support density the tail difference runs from $+1$ to $-1$ as $b$ varies. At the root, $g(b_\kappa,\kappa)=0$ gives
\[
P(Z-b_\kappa>\kappa)-P(Z-b_\kappa<-\kappa)=\frac1\kappa\int_0^\kappa u\,\bigl[p_Z(b_\kappa-u)-p_Z(b_\kappa+u)\bigr]\,du,
\]
so strict skewness is a one-sided condition on the local asymmetry of the density around the Huber location. For asymmetric Laplace noise it can be verified in closed form; for skew-normal (shape $4$) and Gumbel noise we have only numerical evidence (no sign change on $\kappa\in[0.01,40]$, with the difference falling below machine precision at large $\kappa$).

\begin{proposition}[Asymmetric Laplace noise is strictly skewed]
\label{prop:laplace}
Let $Z$ have density $p_Z(z)=Ae^{-rz}$ for $z\ge0$ and $p_Z(z)=Ae^{lz}$ for $z<0$, where $r,l>0$, $r\neq l$ and $A=rl/(r+l)$, or let $Z$ be any translate of such a variable. Then $Z$ is strictly skewed, with $\sigma=-1$ if $r>l$ and $\sigma=+1$ if $r<l$.
\end{proposition}

\begin{proof}
The density is continuous and strictly positive, and $\E|Z|<\infty$. Write $N(b,\kappa):=P(Z-b>\kappa)-P(Z-b<-\kappa)$ and
$T(b,\kappa):=\E[(Z-b)\1\{|Z-b|\le\kappa\}]$, so that $g=T+\kappa N$. Translating $Z$ translates $b_\kappa$ by the same amount and leaves $N(b_\kappa,\kappa)$ unchanged. Replacing $Z$ by $-Z$ swaps $r$ and $l$, replaces $b_\kappa$ by $-b_\kappa$ and changes the sign of $N(b_\kappa,\kappa)$. It therefore suffices to treat the untranslated law with $r>l$ and to show $N(b_\kappa,\kappa)<0$ for every $\kappa>0$. Fix $\kappa$ and write $c:=-b_\kappa$.

\emph{Step 1: $c>0$.} By oddness of $\psi_\kappa$,
$g(0,\kappa)=A\int_0^\infty\psi_\kappa(z)(e^{-rz}-e^{-lz})\,dz<0$, because $\psi_\kappa>0$ on $(0,\infty)$ and $r>l$. Since $g(\cdot,\kappa)$ is strictly decreasing and vanishes at $b_\kappa$, we get $b_\kappa<0$.

\emph{Step 2: $c\ge\kappa$.} Then $[b_\kappa-\kappa,b_\kappa+\kappa]\subseteq(-\infty,0]$,
where $p_Z(z)=Ae^{lz}$, so
\[
T(b_\kappa,\kappa)=Ae^{-lc}\int_{-\kappa}^{\kappa}ue^{lu}\,du
=\frac{2Ae^{-lc}}{l^2}\bigl(l\kappa\cosh(l\kappa)-\sinh(l\kappa)\bigr)>0,
\]
since $x\cosh x>\sinh x$ for $x>0$. As $g(b_\kappa,\kappa)=0$, it follows that
$N(b_\kappa,\kappa)=-T(b_\kappa,\kappa)/\kappa<0$.

\emph{Step 3: $c<\kappa$.} For any $c'\in(0,\kappa)$ we have
$-c'+\kappa>0>-c'-\kappa$, so the two tails of $Z+c'$ beyond $\pm\kappa$ lie in
different exponential pieces and
\[
N(-c',\kappa)=\frac{A}{r}e^{-r(\kappa-c')}-\frac{A}{l}e^{-l(\kappa+c')},
\]
which is negative if and only if $c'<c^*(\kappa):=\bigl((r-l)\kappa+\log(r/l)\bigr)/(r+l)$.
Here $c^*(\kappa)>0$, and $c^*(\kappa)\ge\kappa$ if and only if
$\kappa\le\kappa_0:=\log(r/l)/(2l)$. If $\kappa\le\kappa_0$, then
$c<\kappa\le c^*(\kappa)$ and we are done. Let $\kappa>\kappa_0$, so that
$0<c^*<\kappa$ with $c^*=c^*(\kappa)$, and let
$E_0:=\frac{A}{r}e^{-r(\kappa-c^*)}=\frac{A}{l}e^{-l(\kappa+c^*)}$ be the common value of the two tail probabilities, so that $N(-c^*,\kappa)=0$. Using $\E Z=1/r-1/l$ and integrating the two exponential tails,
\begin{gather*}
\E[Z+c^*]=c^*+\tfrac1r-\tfrac1l,\qquad
\E[(Z+c^*)\1\{Z+c^*>\kappa\}]=E_0\bigl(\kappa+\tfrac1r\bigr),\\
\E[(Z+c^*)\1\{Z+c^*<-\kappa\}]=-E_0\bigl(\kappa+\tfrac1l\bigr),
\end{gather*}
hence
\[
g(-c^*,\kappa)=T(-c^*,\kappa)=c^*-\Bigl(\frac1l-\frac1r\Bigr)(1-E_0)=:\varphi(\kappa).
\]
On $[\kappa_0,\infty)$ the quantity $\kappa-c^*=(2l\kappa-\log(r/l))/(r+l)$ is increasing and vanishes at $\kappa_0$. So $E_0$ decreases from $A/r=l/(r+l)<\tfrac12$,
and
\[
\varphi'(\kappa)=\frac{r-l}{r+l}-\Bigl(\frac1l-\frac1r\Bigr)\frac{2rl}{r+l}E_0
=\frac{r-l}{r+l}\,(1-2E_0)>0.
\]
At $\kappa_0$, where $c^*=\kappa_0$ and $E_0=l/(r+l)$,
\[
\varphi(\kappa_0)=\frac1{2l}\Bigl(\log x-\frac{2(x-1)}{x+1}\Bigr),\qquad x:=\frac rl>1,
\]
which is positive because $\frac{d}{dx}\bigl(\log x-\frac{2(x-1)}{x+1}\bigr)=\frac{(x-1)^2}{x(x+1)^2}>0$
and the bracket vanishes at $x=1$. Hence $g(-c^*,\kappa)=\varphi(\kappa)>0$. As $g(\cdot,\kappa)$ is strictly decreasing, $b_\kappa>-c^*$, that is, $c<c^*(\kappa)$, and therefore $N(b_\kappa,\kappa)<0$. \qed
\end{proof}

With~\eqref{eq:dbkappa} below, it follows that the optimal offset $b_\kappa$ is strictly decreasing in $\kappa$ when the right tail is the lighter one ($r>l$), and that Theorem~\ref{thm:huber} holds without further conditions on every asymmetric Laplace location model.

\begin{theorem}[Canonical Huber thresholds]
\label{thm:huber}
Let $X=\theta+Z$ with $Z$ strictly skewed, and let $\cO_{\mathrm{eq}}$ be oracle
optimality over $\cD_{\mathrm{eq}}$. Then $\cO_{\mathrm{eq}}(L_\kappa)=\{\delta_{b_\kappa}\}$
for every $\kappa>0$, and $\kappa\mapsto b_\kappa$ is $C^1$ and strictly monotone with
\begin{equation}
\label{eq:dbkappa}
\frac{db_\kappa}{d\kappa}=\frac{P(Z-b_\kappa>\kappa)-P(Z-b_\kappa<-\kappa)}{P(|Z-b_\kappa|\le\kappa)}.
\end{equation}
Consequently the rays $\{cL_\kappa:c>0\}$, $\kappa>0$, form a realizable total exclusivity partition of the canonical cone $\mathcal{K}_{\mathrm H}:=\{cL_\kappa:c,\kappa>0\}$.
The theorem concerns only this canonical family.
\end{theorem}

\begin{proof}
\emph{Reduction.} $R_\kappa(b):=\E\,\ell_\kappa(Z-b)$ is convex and finite. Since $\ell_\kappa$ is $\kappa$-Lipschitz, its difference quotients are bounded by $\kappa$, and dominated convergence gives $R'_\kappa(b)=-g(b,\kappa)$. Thus $b$ is optimal if and only if $g(b,\kappa)=0$, that is, $b=b_\kappa$. Uniqueness holds because $\partial_bg<0$ by positivity of the density, and existence follows from the limits $g(b,\kappa)\to\mp\kappa$ as $b\to\pm\infty$ (dominated convergence with
bound $\kappa$).

\emph{Smoothness.} Using
$\psi_\kappa(u)=u\1\{|u|\le\kappa\}+\kappa\,\mathrm{sign}(u)\1\{|u|>\kappa\}$,
\[
g(b,\kappa)=\int_{b-\kappa}^{b+\kappa}(z-b)\,p_Z(z)\,dz+\kappa P(Z-b>\kappa)-\kappa P(Z-b<-\kappa).
\]
Since $p_Z$ is continuous, the Leibniz rule applies. The boundary terms $\pm\kappa\,p_Z(b\pm\kappa)$ cancel against those of the tail terms, and
\[
\partial_bg=-P(|Z-b|\le\kappa),\qquad
\partial_\kappa g=P(Z-b>\kappa)-P(Z-b<-\kappa),
\]
both continuous in $(b,\kappa)$. So $g\in C^1$, and since $\partial_bg<0$ the implicit function theorem gives~\eqref{eq:dbkappa}.

\emph{Injectivity and partition.} By Definition~\ref{def:skewed} the numerator
of~\eqref{eq:dbkappa} has constant sign $\sigma$, so $\kappa\mapsto b_\kappa$ is strictly monotone. Since $\cO_{\mathrm{eq}}(cL_\kappa)=\cO_{\mathrm{eq}}(L_\kappa)$
for $c>0$ (Proposition~\ref{prop:invariance}(i) below), and since $cL_\kappa$ determines $\kappa$ as its maximal quadratic-core radius, Lemma~\ref{lem:invariant} applies with $I(cL_\kappa)=\kappa$. \qed
\end{proof}

\paragraph{Scope and sharpness.}
Theorem~\ref{thm:huber} does not establish exclusivity for the full class of convex losses $\ell(a-\theta)$ whose derivative is linear exactly on $[-\kappa,\kappa]$. Beyond the core, the shape of $\ell'$ changes the root of the estimating equation: one can modify $\ell'$ beyond the core so that $\E\,\ell'(Z-b_{\kappa'})=0$ for a prescribed $\kappa'\neq\kappa$. A loss with core radius $\kappa$ tuned in this way shares its optimal offset with a canonical Huber loss of a different radius, so the core radius alone cannot be expected to be an invariant of that larger class. Some asymmetry is necessary, since under symmetry the whole canonical family collapses; Proposition~\ref{prop:huber-collapse} gives then a second collapse mechanism.

\paragraph{Magnitude.}
In order to see the size of the effect, let $Z$ be a mean-centred asymmetric Laplace variable with density proportional to $e^{-\beta z}$ for $z\ge0$ and $e^{z/\beta}$ for $z<0$, with $\beta=2.5$, before centring. This is the case $r=\beta$, $l=1/\beta$ of Proposition~\ref{prop:laplace}, so $Z$ is strictly skewed with $\sigma=-1$. Table~\ref{tab:huber} gives the roots $b_\kappa$, computed by quadrature and bisection. As $\kappa$ grows, the offset moves from near the noise median ($0.738$) to the mean ($0$). The spread exceeds $0.6$, about $0.3$ noise standard deviations ($\mathrm{sd}(Z)\approx2.53$). The threshold $\kappa$ thus determines the estimand as well as the sensitivity to outliers.

\begin{table}[t]
\centering
\caption{Roots $b_\kappa$ of~\eqref{eq:huber-foc}, i.e.\ the optimal equivariant estimators $X-b_\kappa$, under mean-centred asymmetric Laplace noise
($\beta=2.5$).}
\label{tab:huber}
\begin{tabular}{lccccccc}
\toprule
$\kappa$ & $0.25$ & $0.5$ & $1$ & $2$ & $3$ & $5$ & $\infty$\\
\midrule
$b_\kappa$ & $0.734$ & $0.722$ & $0.672$ & $0.487$ & $0.315$ & $0.133$ & $0$\\
\bottomrule
\end{tabular}
\end{table}

\paragraph{Asymmetry and robustness together.}
M-quantiles \cite{BrecklingChambers1988} combine the two mechanisms of this section through the asymmetric Huber derivative
$\psi_{\kappa,\tau}(u)=2\psi_\kappa(u)\,(\tau\1\{u>0\}+(1-\tau)\1\{u\le0\})$. This family contains quantiles ($\kappa\to0$) and expectiles ($\kappa\to\infty$)
\cite{NeweyPowell1987}. One may ask whether
$(\tau,\kappa)$ is jointly an invariant; we return to this question in Section~\ref{sec:discussion}.

\section{Structural properties of the loss space}
\label{sec:structure}

The three families have features in common that the individual examples do not account for. In each, rescaling a loss never changes its optimal rules, while mixing two losses with different optimal rules can produce a loss with a new optimal rule. Both facts are proved here in general, together with a description of the effect of changing the operator.

\subsection{Compatibility components}
\label{subsec:components}

Fix $\cM\subseteq\cL$ and let $G_{\cO}(\cM)$ be the graph on $\cM$ with an edge between $L$ and $L'$ whenever they are $\cO$-compatible. Its connected components are the \emph{compatibility components} of $\cM$.

\begin{proposition}[Exclusivity and components]
\label{prop:components}
\begin{enumerate}
\item $\mathcal{C}\subseteq\cM$ is an exclusivity region within $\cM$ if and only if it is a union of compatibility components.
\item The compatibility components form the finest total exclusivity partition of $\cM$; a component fails to be realizable only if it is a single loss $L$ with $\cO(L)=\emptyset$.
\item If $\cO(L)$ is a singleton for every $L\in\cM$, compatibility is an equivalence relation and the components are the fibres of $L\mapsto\cO(L)$.
\end{enumerate}
\end{proposition}

\begin{proof}
(i) $\mathcal{C}$ is a region if and only if no edge joins $\mathcal{C}$ to
$\cM\setminus\mathcal{C}$, which holds if and only if $\mathcal{C}$ is a union of
components. (ii) By (i), every total exclusivity partition consists of unions of
components, so it is refined by the partition into components. A component with an edge contains a loss with $\cO(L)\neq\emptyset$. (iii) If $\cO(L)=\{\delta_L\}$, then $L$ and $L'$ are compatible if and only if
$\delta_L=\delta_{L'}$. \qed
\end{proof}

Without (iii), compatibility doesn't need to be transitive, and 
\textit{class} must mean component. Let $Q=\tfrac12\delta_0+\tfrac12\delta_1$. Then
$\int\rho_\tau(\theta-a)\,Q(d\theta)=\tfrac12[\tau+a(1-2\tau)]$ on $[0,1]$, so
$\cO_Q(\rho_\tau)=\{0\}$ for $\tau<\tfrac12$, $[0,1]$ for $\tau=\tfrac12$, and
$\{1\}$ for $\tau>\tfrac12$. Thus $\rho_{1/4}$ and $\rho_{3/4}$ are incompatible, yet
both are compatible with $\rho_{1/2}$, and all three lie in one component. Two
features of $Q$ act together. The gap of $Q$ on $(0,1)$ makes $\cO_Q(\rho_{1/2})$
set-valued, and the atoms at its ends make the neighbouring levels land on its endpoints. With a continuous distribution function, a gap alone makes acts set-valued but leaves compatibility transitive. Assumption~\ref{ass:posterior}
excludes both, and under it $\cO_\Pi$ is single-valued on $\cL_{\mathrm q}^\Pi$ and the components are the level sets of
$\tau$.

\subsection{Invariances of the optimality map}
\label{subsec:invariance}

\begin{proposition}[Invariances]
\label{prop:invariance}
Let $c,c'>0$ and let $h:\Theta\to[0,\infty)$ be measurable.
\begin{enumerate}
\item If $\E[h(\theta)\mid X]<\infty$ a.s., then $\cO_\Pi(cL+h)=\cO_\Pi(L)$; and
$\cO(cL+h)=\cO(L)$ for oracle optimality over any $\cD_0$.
\item Minimax optimality satisfies $\cO(cL)=\cO(L)$, but in general not
$\cO(L+h)=\cO(L)$.
\item If $L,L'$ are $\cO_\Pi$-compatible and their posterior expected losses are
finite, then $\cO_\Pi(cL+c'L')=\cO_\Pi(L)\cap\cO_\Pi(L')$. In particular, when
$\cO_\Pi$ is single-valued, each fibre $\{L:\cO_\Pi(L)=\{\delta\}\}$ (restricted
to losses with finite posterior expected loss) is a convex cone.
\end{enumerate}
\end{proposition}

\begin{proof}
(i) The posterior expected loss of $cL+h$ is $c$ times that of $L$ plus the
finite constant $\E[h(\theta)\mid X=x]$, and the risk of $cL+h$ at $\theta$ is $cR_L(\theta,\delta)+h(\theta)$; neither change affects the minimizers. (ii)
$\sup_\theta cR=c\sup_\theta R$. For the failure under shifts, take
$\Theta=\{0,1\}$, no data, constant rules, $L=(a-\theta)^2$ and $h(0)=0$,
$h(1)=1$: the minimax act is $\tfrac12$ for $L$ and $1$ for $L+h$. (iii) For finite functions $\varphi,\varphi'$ whose minimizer sets intersect, $\arg\min(c\varphi+c'\varphi')=\arg\min\varphi\cap\arg\min\varphi'$,
since $c\varphi+c'\varphi'\ge c\min\varphi+c'\min\varphi'$ with equality exactly on the intersection. Apply this to the posterior expected losses at $P_X$-a.e.\ $x$; compatibility guarantees a nonempty intersection a.s. \qed
\end{proof}

By (i), the parametrization~\eqref{eq:gpl} already absorbs the standard equivalence $L'=cL+h(\theta)$ into the quantile family: $cL_{\tau,g,h}=L_{\tau,cg,ch}$.
The geometric results below use only scaling, which is the invariance shared by all three operators; shift invariance fails for minimax, and the invariants $\tau$, $H_\phi$ and $\kappa$ are already shift-free.

\begin{proposition}[Completion along rays]
\label{prop:cone}
Let $K\subseteq\cL$ be a cone and $\cO$ invariant under positive scaling. Let $\{\mathcal{C}_i\}_{i\in I}$ be pairwise disjoint exclusivity regions within $K$,
each containing a loss with nonempty $\cO$, such that no element of $\mathcal{C}_i$ is a positive multiple of an element of $\mathcal{C}_j$ for $i\neq j$. Then the saturations
$\widetilde{\mathcal{C}}_i:=\{cL:c>0,\,L\in\mathcal{C}_i\}$ form a realizable total exclusivity partition of $K':=\bigcup_i\widetilde{\mathcal{C}}_i$.
\end{proposition}

\begin{proof}
If $cL=c'L'$ with $L\in\mathcal{C}_i$ and $L'\in\mathcal{C}_j$, then $L'=(c/c')L$, so $i=j$. Hence the saturations are disjoint, and they cover $K'$ by construction. Suppose $\delta\in\cO(cL)\cap\cO(L'')$ with $L\in\mathcal{C}_i$
and $L''\in K'\setminus\widetilde{\mathcal{C}}_i$. Write $L''=c''M$ with $M\in\mathcal{C}_j$, $j\neq i$. Then $\delta\in\cO(L)\cap\cO(M)$ by scale invariance, and $M\in K\setminus\mathcal{C}_i$, therefore contradicting exclusivity of
$\mathcal{C}_i$. Realizability is (trivially) inherited. \qed
\end{proof}

The non-homothety hypothesis holds automatically in Section~\ref{sec:families}, because $\tau(cL)=\tau(L)$, $H_{c\phi}=H_\phi$, and $cL_\kappa$ has core radius
$\kappa$. Exclusivity established on normalized representatives therefore extends along rays.

\subsection{Convex interpolation}
\label{subsec:convexity}

Proposition~\ref{prop:invariance}(iii) says that positive combinations of
compatible losses keep their common optimal rules. For two \emph{incompatible}
losses the situation is different. In the mixture example of
Section~\ref{subsec:quantile}, every interior combination has a new Bayes rule,
and the segment passes through a continuum of compatibility classes.
Theorem~\ref{thm:convexity} shows that this is a general feature of Bayes acts,
and nothing specific to quantiles is involved.

Fix a law $Q$ on $\Theta$, let $\cA$ be a closed subset of a finite-dimensional
normed space, and write
\[
B_L(a):=\int L(\theta,a)\,Q(d\theta),\qquad \cO_Q(L)=\arg\min_{a\in\cA}B_L(a).
\]
In the posterior case the theorem applies at each fixed $x$. If its assumptions
hold for $P_X$-a.e.\ $x$ with $a^*_x(0)\neq a^*_x(1)$ on a set of positive
$P_X$-measure, and $x\mapsto a^*_x(\alpha)$ is measurable (e.g.\ by the measurable
maximum theorem), then $\alpha\mapsto\delta_\alpha:=a^*_\cdot(\alpha)$ is continuous
for convergence in $P_X$-probability. Its image is a connected metric space with
at least two points, hence uncountable.

\begin{assumption}[Regular endpoint risks]
\label{ass:endpoints}
$L_0,L_1\in\cL$, and $L_\alpha:=(1-\alpha)L_0+\alpha L_1$ for $\alpha\in[0,1]$.
\begin{enumerate}
\item $B_{L_0}$ and $B_{L_1}$ are finite, continuous and coercive on $\cA$
(that is, $B_{L_i}(a)\to\infty$ as $\|a\|\to\infty$ in $\cA$).
\item For every $\alpha\in[0,1]$, $B_\alpha:=B_{L_\alpha}$ has a unique minimizer
$a^*(\alpha)$.
\item $a^*(0)\neq a^*(1)$.
\end{enumerate}
\end{assumption}

\begin{theorem}[Convex interpolation meets a continuum of classes]
\label{thm:convexity}
Under Assumption~\ref{ass:endpoints}, $\alpha\mapsto a^*(\alpha)$ is continuous and its
image is uncountable. Consequently $\cO_Q$ is single-valued on the segment
$\mathcal{S}:=\{L_\alpha:\alpha\in[0,1]\}$, the compatibility components of
$\mathcal{S}$ are the fibres $\{L_\alpha:a^*(\alpha)=a\}$, and the compatibility
partition of $\mathcal{S}$ has uncountably many classes. In particular, no
countable set of actions contains an optimal act for every $L_\alpha$.
\end{theorem}

\begin{proof}
Each $B_\alpha=(1-\alpha)B_{L_0}+\alpha B_{L_1}$ is finite and continuous, and since
$B_\alpha\ge\min(B_{L_0},B_{L_1})$ it is coercive. Thus $a^*(\alpha)$ exists, and
by (ii) it is unique. Moreover $(\alpha,a)\mapsto B_\alpha(a)$ is jointly
continuous, being affine in $\alpha$ with continuous coefficients.

\emph{Continuity.} Let $\alpha_n\to\alpha$ and $a_n:=a^*(\alpha_n)$. Then
$\min(B_{L_0},B_{L_1})(a_n)\le B_{\alpha_n}(a_n)\le B_{\alpha_n}(a^*(\alpha))\to B_\alpha(a^*(\alpha))$.
The left-hand function is coercive, so $(a_n)$ is bounded. Let $\bar a$ be any
limit point; it lies in $\cA$ because $\cA$ is closed. For every $a\in\cA$,
$B_{\alpha_n}(a_n)\le B_{\alpha_n}(a)$, and passing to the limit along the
subsequence gives $B_\alpha(\bar a)\le B_\alpha(a)$ by joint continuity. Hence
$\bar a=a^*(\alpha)$ by uniqueness. A bounded sequence in a finite-dimensional
space whose only limit point is $a^*(\alpha)$ converges to it.

\emph{Uncountability.} Choose a linear functional $\ell$ with
$\ell(a^*(0))\neq\ell(a^*(1))$. Then $\ell\circ a^*$ is continuous on $[0,1]$, and
by the intermediate value theorem it takes every value between
$\ell(a^*(0))$ and $\ell(a^*(1))$. So the image of $a^*$ is uncountable.

\emph{Classes.} $\cO_Q(L_\alpha)=\{a^*(\alpha)\}$, so by
Proposition~\ref{prop:components}(iii) the components of $\mathcal{S}$ are the
fibres of $a^*$. These are indexed by the uncountable image, and any set of
actions containing an optimal act for each $L_\alpha$ contains that image. \qed
\end{proof}

\paragraph{What the theorem rules out, and what it does not.}
It rules out a finite (or countable) menu of procedures that serves every member
of a convex family of losses satisfying Assumption~\ref{ass:endpoints}. In
this sense the compatibility structure is \emph{conic, not convex}. Rays and
positive combinations within a class stay in the class
(Proposition~\ref{prop:invariance}), but a segment between classes meets
uncountably many of them. The theorem does not say that $\alpha\mapsto a^*(\alpha)$
is injective. The path may remain in an endpoint class on a subinterval (for
$B_0(a)=|a|$ and $B_1(a)=(a-1)^2$ one has $a^*(\alpha)=0$ for
$\alpha\le\tfrac13$), and in general it may return to a fibre it has left. For
the weighted pinball objective of Section~\ref{subsec:quantile} the path is
injective, and it leaves both endpoint classes for every $\alpha\in(0,1)$. Under Assumption~\ref{ass:posterior}
with $\int|\theta|\,dQ<\infty$ and $\tau_0<\tau_1$, the segment from $\rho_{\tau_0}$ to $\rho_{\tau_1}$
satisfies Assumption~\ref{ass:endpoints}, and its Bayes acts sweep the interval of
quantiles $[q(\tau_0),q(\tau_1)]$.

The theorem does not say that convex families lack finite exclusivity
\emph{regions}. Grouping the $L_\alpha$ by finitely many Borel sets of act values
gives a finite region partition (Proposition~\ref{prop:components}(i)), each
of whose parts is a union of possibly uncountably many classes. The theorem also
needs uniqueness. When Bayes
acts are set-valued, as for the median of a distribution with a gap, the
argument has to be redone with set-valued continuity, and the conclusion
concerns components; without uniqueness the conclusion can fail outright (for
$Q=\delta_0$, $L_0=|\theta-a|$ and $L_1=|\theta+1-a|$ the whole segment is one
component). Finally, the theorem says nothing about minimax or
restricted-oracle operators, for which continuity of the optimal act in $\alpha$
must be established separately.

\subsection{Coarsening the operator}
\label{subsec:coarsening}

Exclusivity is relative to the operator, and statistical practice often reads
only part of an optimal rule: the sign of a score, a derived probability, or a
summary of an estimator. We write $\cO_1\preceq\cO_2$ if
$\cO_1(L)\subseteq\cO_2(L)$ for all $L$. A \emph{readout} is a fixed map
$\pi:\cD\to\cD'$, the same for every loss, with induced operator
$(\pi\circ\cO)(L):=\pi(\cO(L))$. A loss-indexed decoding $L\mapsto\pi_L$ is not a
readout in this sense and is not covered by
Proposition~\ref{prop:coarsening}.

\begin{proposition}[Coarsening can only merge components]
\label{prop:coarsening}
Let $\cM\subseteq\cL$.
\begin{enumerate}
\item Let $\cO_1\preceq\cO_2$. Compatibility under $\cO_1$ implies
compatibility under $\cO_2$. Every $\cO_1$-component of $\cM$ lies in an $\cO_2$-component,
and every exclusivity region within $\cM$ under $\cO_2$ is one under $\cO_1$.
\item The same conclusions hold with $\cO_1=\cO$ and $\cO_2=\pi\circ\cO$, for any
readout $\pi$.
\end{enumerate}
\end{proposition}

\begin{proof}
(i) If $\delta\in\cO_1(L)\cap\cO_1(L')$, then
$\delta\in\cO_2(L)\cap\cO_2(L')$. So $G_{\cO_1}(\cM)$ is a subgraph of
$G_{\cO_2}(\cM)$, and adding edges can only merge components. A set with no
$\cO_2$-edge leaving it has no $\cO_1$-edge leaving it. (ii) If
$\delta\in\cO(L)\cap\cO(L')$, then $\pi(\delta)\in\pi(\cO(L))\cap\pi(\cO(L'))$,
and the rest is as in (i). \qed
\end{proof}

\section{Collapse and operator dependence}
\label{sec:collapse}

The exclusivity results of Section~\ref{sec:families} hold for a particular operator
and a particular class of models. When either is changed they can fail in two
ways: the operator may read less of the optimal rule, or the model may give the
loss nothing to distinguish.

\subsection{Score, probability and classifier}
\label{subsec:margin-collapse}

For a margin loss, an optimal score can be used in three ways: as a raw score,
as a classifier via the fixed readout $\mathrm{sign}$, or as a probability via the
loss-indexed decoding map $d_\phi(t):=H_\phi(t)/(1+H_\phi(t))$, which depends on the
loss through its link. Let
$\cL_{\mathrm{cal}}$ be the convex classification-calibrated margin losses
\cite{BartlettJordanMcAuliffe2006} (hinge, logistic, exponential, squared,
\ldots). Let $\cO_{\mathrm{score}}$ be as in Theorem~\ref{thm:margin}, extended to
$\cL_{\mathrm{cal}}$, and let $\cO_{\mathrm{class}}:=\mathrm{sign}\circ\cO_{\mathrm{score}}$,
with sign rules identified when they agree $P_X$-a.s.\ off the tie set
$T=\{\eta=\tfrac12\}$. For $L_\phi\in\cL_{\mathrm m}$, let
$\cO_{\mathrm{prob}}(L_\phi):=\{d_\phi\circ f:f\in\cO_{\mathrm{score}}(L_\phi)\}$.

\begin{proposition}[Score, probability and classifier]
\label{prop:readouts}
Assume $\eta(X)\in(0,1)$ $P_X$-a.s.
\begin{enumerate}
\item Under Assumption~\ref{ass:richness}, $\cL_{\mathrm m}$ is partitioned by
$H_\phi$ under $\cO_{\mathrm{score}}$ (Theorem~\ref{thm:margin}).
\item $\cO_{\mathrm{prob}}(L_\phi)=\{\eta\}$ for every $L_\phi\in\cL_{\mathrm m}$, so
$\cL_{\mathrm m}$ is a single class under $\cO_{\mathrm{prob}}$.
\item $\cO_{\mathrm{class}}(L_\phi)=\{\mathrm{sign}(\eta-\tfrac12)\}$ for every
$L_\phi\in\cL_{\mathrm{cal}}$, so $\cL_{\mathrm{cal}}$ is a single class under
$\cO_{\mathrm{class}}$.
\end{enumerate}
\end{proposition}

\begin{proof}
(ii) By~\eqref{eq:margin-foc}, $H_\phi(f^*_\phi)=\eta/(1-\eta)$. (iii) For convex
calibrated $\phi$ one has $\phi'(0)<0$ \cite{BartlettJordanMcAuliffe2006}, so
$\phi(t)\to\infty$ as $t\to-\infty$. For $\eta\in(0,1)$, $C_\eta$ is then
convex, continuous and coercive, so its minimizer set is a nonempty compact
interval depending upper hemicontinuously on $\eta$, and a measurable minimizing
score rule exists by the measurable selection theorem. By calibration every
minimizer has sign $\mathrm{sign}(\eta-\tfrac12)$ off $T$. \qed
\end{proof}

Part (iii) is the known fact that calibrated surrogates share the Bayes
classifier \cite{Zhang2004,BartlettJordanMcAuliffe2006}, and we do not present
it as new. The framework places it in a general setting: part (iii) is an
instance of
Proposition~\ref{prop:coarsening}(ii) with $\pi=\mathrm{sign}$: the
score-level classes of Theorem~\ref{thm:margin} can only merge under the sign
readout, and by (iii) they merge into a single class. Part (ii) is \emph{not} an
instance of Proposition~\ref{prop:coarsening}. The decoding $d_\phi$ depends on
the loss, and the collapse there is a direct consequence
of~\eqref{eq:margin-foc}: each loss decodes its own Bayes score back to $\eta$.
The logistic family of Section~\ref{subsec:margin} exhibits all three
behaviours. Whether the choice of margin loss matters is therefore not a
property of the losses themselves; it depends on what the downstream task reads,
namely the raw score on a fixed scale, a probability, or a label.

\subsection{Huber thresholds without asymmetry}
\label{subsec:huber-collapse}

Theorem~\ref{thm:huber} needs strict skewness, and collapses occur without it. Under
symmetric noise all canonical thresholds share one estimator
(Section~\ref{subsec:huber}). In the posterior Bayes setting a second mechanism
appears.

\begin{proposition}[Bounded posterior residuals]
\label{prop:huber-collapse}
Let $\bar\kappa>0$, and suppose that for $P_X$-a.e.\ $x$ the posterior $\Pi_x$ is
supported on a closed interval $I_x$ of length at most $\bar\kappa$, and that
$\cD$ contains the posterior-mean rule $x\mapsto\E[\theta\mid X=x]$. Then
$\cO_\Pi(L_\kappa)=\cO_\Pi(L_{\kappa'})$ for all $\kappa,\kappa'\ge\bar\kappa$, so
$\{L_\kappa:\kappa\ge\bar\kappa\}$ lies in a single class.
\end{proposition}

\begin{proof}
Fix such an $x$. Since $\ell_\kappa$ is strictly increasing in $|u|$, any
$a>\max I_x$ has strictly larger posterior risk than $\max I_x$, and similarly
below. Hence all minimizers lie in $I_x$. For $a,\theta\in I_x$ we have
$|a-\theta|\le\bar\kappa\le\kappa$, so $\ell_\kappa(a-\theta)=\tfrac12(a-\theta)^2$.
The posterior risks of $L_\kappa$ and $L_{\kappa'}$ therefore coincide on $I_x$,
and both are minimized uniquely at the posterior mean, so
$\cO_\Pi(L_\kappa)=\{\E[\theta\mid X]\}$ for every $\kappa\ge\bar\kappa$. \qed
\end{proof}

Here the threshold has no effect because the loss never leaves its quadratic
part on the relevant support. In both collapses identification fails, in the
sense that the Bayes act does not determine the loss parameter. Theorem~\ref{thm:huber}
identifies the model class (strictly skewed, full-support noise) on which it
does.

\section{Discussion}
\label{sec:discussion}

Relative to a model, a decision class and an optimality operator, the paper
establishes the following. Exclusivity regions are exactly unions of
compatibility components. For three standard families, a single loss-side
quantity determines the component: the quantile level for strictly consistent
quantile scores under posterior Bayes optimality and Assumption~\ref{ass:posterior};
the margin invariant (link) for regular margin losses at the score level under
Assumption~\ref{ass:richness}; and the threshold for the canonical Huber family
among equivariant rules on strictly skewed location models, which include all
asymmetric Laplace models. Optimal rules are unchanged by positive scaling and
by positive combination of compatible losses. A segment between losses with
distinct Bayes acts meets uncountably many compatibility classes when the endpoint
risks are coercive and every Bayes act along it is unique. Coarser operators have
coarser components.

This viewpoint is complementary to that of elicitation, which fixes a functional
and asks which losses target it. Here a loss family is fixed and we ask how the
Bayes-act map partitions it. Where the family consists of consistent scores, as
for quantiles, the partition is the elicitation result read backwards, and
Theorem~\ref{thm:quantile} states in addition that the level is the
\emph{only} feature of the score that the Bayes rule reads. Where the family does
not target a single functional, as for margin losses read as scores,
probabilities or labels, the partition records information that consistency
alone does not.

The results have several limitations. All of them are model-relative: the
quantile and margin
theorems need posterior or class-probability regularity, and the two-point
example of Section~\ref{subsec:components} and the constant-$\eta$ example after
Theorem~\ref{thm:margin} show that the conclusions fail when these are weakened.
Set-valued optima
destroy transitivity, and then only components are canonical
(Section~\ref{subsec:components}). The Huber result covers the canonical family
only; we do not expect the quadratic-core radius to be an invariant of the larger
class of convex losses with that radius. Strict skewness, finally, is proved for
asymmetric Laplace noise (Proposition~\ref{prop:laplace}) but only checked
numerically for the skew-normal and Gumbel laws.

We mention two open directions. The first is \emph{model-robust} exclusivity:
partitions that persist uniformly over a model class, such as all location models
with strictly skewed full-support noise. For quantiles this holds by
Theorem~\ref{thm:quantile}. For Huber thresholds, Proposition~\ref{prop:laplace}
gives one such class (all asymmetric Laplace location models), but a
characterization of the noise laws on which it holds is open. The second is
multidimensional families, of which M-quantile losses
\cite{BrecklingChambers1988} are the natural first case. Is $(\tau,\kappa)$ a
complete invariant of the asymmetric Huber family, and on which models? A positive
answer would combine the quantile and Huber results. A negative answer would
exhibit a curve in parameter space along which asymmetry and robustness trade off
without changing the estimand.

\bibliographystyle{spmpsci}
\bibliography{bibliography-metrika}

\end{document}